\documentclass[11pt,a4paper,twoside]{amsart}

\usepackage{amsfonts, amssymb,indentfirst}
\usepackage[all]{xy}

\newtheorem{prop}{Proposition}[section]
\newtheorem{theorem}[prop]{Theorem}
\newtheorem{lemma}[prop]{Lemma}

\newtheorem{coro}[prop]{Corollary}
\newtheorem{remark}[prop]{Remark}

\newenvironment{rem}{\begin{remark}\rm}{\end{remark}}
\newcommand{\cqd}{\hfill$\Box$}
\newcommand{\mor}[0]{\operatorname{Mor}}

\title[ON BRILL-NOETHER LOCI OVER QUOT SCHEMES AND A TORELLI THEOREM]
{ON BRILL-NOETHER LOCI OVER QUOT SCHEMES AND A TORELLI THEOREM}
\author[Cristina Mart{\'\i}nez]{Cristina Mart{\'\i}nez}

\date{October 2011}
\begin{document}
\maketitle

\begin{abstract}

 We prove a  non abelian Torelli type result for  smooth projective curves by working
 in the derived category of some associated polarized Quot schemes and defining Brill-Noether loci and Abel-Jacobi maps on it.

\end{abstract}

\section{Introduction}




The derived category of coherent sheaves on a variety $X$ defined over an algebraically closed field $k$, is a triangulated category whose objects
are bounded and coherent cohomology sheaves on $X$. Due to a result of Orlov (\cite{Orl}), an equivalence
$F:\mathcal{D}(X)\rightarrow \mathcal{D}(X')$ between derived categories of coherent
sheaves on smooth projective varieties $X, X'$ is always of
Fourier-Mukai type, that is, there exists a unique (up to
isomorphism) object $\mathcal{P}\in \mathcal{D}(X\times X')$ such that the
functor $F$ is isomorphic to the derived functor of the pushforward map:
$$\Xi_{\mathcal{P}}(-):=\mathbf{ R}q_{*}(\mathcal{P}\otimes p^{*}(-)),$$ where $p$ and $q$
are the projections of $X\times X'$ onto $X$ and $X'$
respectively.
For smooth projective curves, a derived equivalence always
corresponds to an isomorphism. In particular this implies the classical Torelli
Theorem. 
If there is an equivalence between the derived categories of two
smooth projective curves, then there is an isomorphism between the
Jacobians of the curves that preserves the principal polarisation,
\cite{Be}.



Let $C$ be a projective irreducible non-singular curve over an
algebraic closed field $k$. 
\noindent We consider the Jacobian $J(C)$ of $C$ which is a projective abelian
variety parametrizing line bundles of degree 0 on $C$. The Brill-Noether locus as defined in \cite{ACGH}  is the subvariety of $J(C)$ parametrizing
varieties of special divisors. The classical Torelli Theorem allows
us to recover the curve from its Jacobian as a polarized abelian
variety.


In general, it is of interest to know how much information about a
space $X$ can be recovered from Hodge data on $X$. In the language of motives in the sense of D. Arapura in \cite{Ar}, the classical Torelli theorem implies that the Jacobian of a curve $C$ is motivated by the curve, that is, the motive of $C$ is contained in the category generated from $J(C)$ by taking sums, summands and products. This generalizes to other moduli spaces. For the moduli space of vector bundles over a curve, this result was first proved by
S. del Ba{\~n}o  in
\cite{Ba}. 
 
\noindent The space of matrix divisors (of given rank and degree) studied by Bifet-Ghione-Letizia in \cite{BGL}, is closely related to the moduli space of stable rank $n$ and degree $d$ bundles on $C$, which since work of Weil has been considered a non-abelian variant of the Jacobian.

\noindent Here we show that a
smooth projective curve $C$ over $k$ is determined by a certain
Quot scheme compactification of the scheme of degree $d$ morphisms
from the curve to the Grassmannian $G(2,4)$ with a certain
polarization on this Quot scheme. The proof uses a Fourier-Mukai
functor  along the lines of the Beilinson-Polishchuk proof of the
classical Torelli theorem, \cite{BP}. In the genus 0 case, these
spaces are considered in \cite{Mar1} as parameter spaces for
rational ruled surfaces in order to solve a certain enumerative
problem. However the Fourier-Mukai functor is defined on a
general Quot scheme parametrizing quotient sheaves of a trivial
bundle on $C$. This method was also applied to Prym varieties by
J. C. Naranjo in \cite{Nar}.

\subsubsection*{Conventions} 
Points of a scheme are always closed points. By a sheaf on a scheme $X$ we mean a coherent $\mathcal{O}_{X}-$module.
By $\mathcal{D}(X)$ we denote the derived category of cohomology sheaves on a smooth projective variety $X$ over an algebraically closed field $k$, (the case of primary interest for us is when $k$ is of characteristic 0), and by $\mathcal{D}^{b}(X)$ its full subcategory of bounded complexes. We write $H^{i}(X)$ for its cohomology  with rational coefficients. 
For $E\in
\mathcal{D}(X)$ a rank $n$ vector bundle on $X$, we denote by $det\,(E)$ its
determinant line bundle.

\section{Geometry of the Jacobian and Quot schemes}
\subsection{Non abelian Torelli problem for smooth projective curves.} 
We fix  a complete non-singular genus $g$ curve $C$ over  an
algebraically closed field $k$. Let $Pic\,(C)$ be
 the Picard group of $C$ parametrizing all line bundles  over $C$, and  $Pic^{d}(C)$  the degree $d$ subset in it. 
The Jacobian $J=J(C)$ of $C$ is an abelian variety such that
the group of $k-$points of $J$ is isomorphic to $Pic^{0}(C)$ (resp.
$J^{d}(C)$ is isomorphic to $Pic^{d}(C)$).  
There is a canonical map $C\rightarrow Pic^{1}(C)$ (see \S I. 3.3.1 of \cite{ACGH}), and upon choice of a point $P\in C$, we get a map $e: C\rightarrow J(C)$, normalized such that $e(P)=0$. 

For every $d>0$, we denote by $Sym^{d}C$ the $d$th symmetric power of a curve
$C$. By definition, $Sym^{d}C$ is the quotient of $C^{d}$ by the action of the
symmetric group $S_{d}$. We can identify the set of effective divisors of degree $d$ on $C$
 with the set of $k-$rational points of the symmetric
power $Sym^{d}C$, that is, $Sym^{d}C$ represents the functor of
families of effective divisors of degree $d$ on $C$. 

\begin{theorem}(Torelli) Let $C_{1}$ and $C_{2}$ be two smooth projective curves of
genus $g>1$ over $k$. If there is an isomorphism between the
Jacobians $J(C_{1})$ and $J(C_{2})$ preserving the principal
polarization then $C_{1}\cong C_{2}$.
\end{theorem}

The subset in $Pic^{g}(C)$ consisting of line bundles $L$ with
$h^{0}(L)=1$ corresponds to the set of $k-$points of an open
subset in $Sym^{g}C$. Translating this subset by various line
bundles of degree $-g$ we obtain algebraic charts for
$Pic^{0}(C)$. The Jacobian variety $J$ is constructed by gluing
together these open charts. It is a consequence of Torelli's
theorem that if
$$Sym^{g-1}C_{1}\cong Sym^{g-1}C_{2}, \ \ then \ C_{1}\cong
C_{2}.$$

The next theorem states that the same result continues to hold for
all $d\geq 1$ with one exception:

\begin{theorem}\cite{Fak} Let $C_{1}$ and $C_{2}$ be two smooth projective curves of genus
$g\geq 2$ over an algebraically closed field $k$. If
$Sym^{d}C_{1}\cong Sym^{d}C_{2}$ for some $d\geq 1,$ then
$C_{1}\cong C_{2}$ unless $g=d=2$.
\end{theorem}
It is well known that there exist non-isomorphic curves of genus 2 with isomorphic Jacobians. Everett W. Howe provides some examples in \cite{Ho}.

\subsection{Varieties of special divisors} 
There is a fundamental relation between linear series on $C$ and maps of $C$ to projective spaces expressed in the language of line bundles. Given a non-zero holomorphic section $s$ of  a line bundle $L$, and denoting by $D$ the divisor of $s$, there is an isormophism between $L$ and $\mathcal{O}(D)$. The complete linear series $\mathcal{L}=\mathbb{P}H^{0}(C,L)$ is the set of effective divisors equivalent to $D$. A linear series $\mathbb{P}V$, where $V$ is a vector space of $\mathcal{L}$ is said to be a $g^{r}_{d}$ if $deg\,(D)=d$ and $dim\,(V)=r+1$. 

There exists an algebraic variety $G^{r}_{d}$ parametrizing the series $g^{r}_{d}$. 
We observe that as $g=0$, the space of $g^{r}_{d}$ is precisely the Grassmannian $G(r,d)$ of $r-$planes in a $d-$dimensional vector space. This allows us to see $G^{r}_{d}$ as a generalized Grassmannian, (see \S IV. 3 of \cite{ACGH}).

The variety of special divisors 
$W^{r}_{d}\subset J^{d}(C)=J(C)$ parametrizes line bundles $L$ of
degree $d$ such that $h^{0}(L)>r$.

One has a canonical scheme
structure on $W^{r}_{d}$, since it can be described as the
degeneration locus of some morphism of vector bundles on $J^{d}$, (see \S IV. 1 of \cite{ACGH}).  
The subscheme $W^{0}_{g-1}$ is exactly the theta divisor
$\Theta\subset J^{g-1}$. All theta divisors in the Jacobian are
translations of the natural divisor $\Theta \subset J^{g-1}$. We
have a canonical involution corresponding to the map $\nu:
\Theta\rightarrow \Theta,$ $L \rightarrow K_{C}\otimes L^{-1}$,
where $K_{C}$ denotes the canonical line bundle over the curve $C$
and $L^{-1}$ is the dual line bundle of $L$. There is a canonical
identification of $Pic^{0}(J^{g-1})$ with $Pic^{0}(J)=\widehat{J}$
induced by any standard isomorphism $J\rightarrow J^{g-1}$ given
by some line bundle of degree $g-1$.
 The corresponding Fourier transform $\mathcal{F}$ on the
derived categories of coherent sheaves on $\widehat{J}$ and
$J^{g-1}$ is an equivalence.

Denote by $\Theta^{ns}$ the open subset of smooth points of
$\Theta$. We can identify $\Theta^{ns}$ with an open subset of
$Sym^{g-1}C$ consisting of effective divisors $D$ of degree $g-1$,
such that $h^{0}(D)=1$. 

For sufficiently large degree $d$, the morphism $\sigma^{d}$
$$\begin{array}{lll}Sym^{d}C & \stackrel{\sigma^{d}}{\rightarrow} & J^{d} \\
D & \rightarrow & {\rm isomorphism\ class \ of\ }
\mathcal{O}_{C}(D)
\end{array}$$
is a projective bundle. The fiber of $\sigma^{d}$ over $L$ is the
variety of effective divisors $D$ such that
$\mathcal{O}_{C}(D)\cong L$. Further, $(\sigma^{d})^{-1}(L)\cong
\mathbb{P}H^{0}(C,L)$. Let us identify $J^{d}$ with $J$ using the
line bundle $\mathcal{O}_{C}(dp)$, where $p\in C$ is a fixed
point. Then we can consider $\sigma^{d}$ as a morphism
$Sym^{d}C\rightarrow J$ sending $D$ to $\mathcal{O}_{C}(D-dp)$. In
more invariant terms: Let $L$ be a line bundle of degree $d>2g-2$
on $C$. Then the morphism $\sigma_{L}: Sym^{d}C\rightarrow J$
sending $D$ to $\mathcal{O}_{C}(D)\otimes L^{-1}$ can be
identified with the projective bundle associated with
$\mathcal{F}(L)$, the Fourier transform of $L$. 


\subsection{Higher rank divisors} Let $\mathcal{O}_{C}$ be the structure sheaf of the curve $C$
and let $K$ be its field of rational functions, considered as a
constant $\mathcal{O}_{C}$ module. Following \cite{BGL}, we define a
divisor of rank $r$ and degree $d$ or $(r,d)$ divisor as a coherent
sub $\mathcal{O}_{C}$-module of $K^{r}= K^{\oplus r}$, having rank
$r$ and degree $d$.

This set can be identified with the set of rational points of an
algebraic variety $Div^{r,d}_{C/k}$ which may be described as
follows. For any effective ordinary divisor $D$, set:
$$Div^{r,d}_{C/k}(D)=\{E\in Div^{r,d}_{C,k}|E\subset \mathcal{O}_{C}(D)^{r}\},$$
where $\mathcal{O}_{C}(D)$ is considered as a submodule of $K^{r}$.

The space of all matrix divisors of rank $r$ and degree $d$  can be identified with the set of
rational points of $Quot^{m}_{\mathcal{O}_{C}(D)^{r}/C/k}$
parametrizing torsion quotients of $\mathcal{O}_{C}(D)^{r}$ and
having degree $m=r\cdot deg\,D-d$. It is a smooth projective
irreducible variety. As in the Jacobian case, tensoring by
$\mathcal{O}_{C}(-D)$ defines an isomorphism between
$\mathcal{Q}_{r,d}(D)=Quot^{m}_{\mathcal{O}_{C}(D)^{r}/C/k}$ and
$Quot^{m}_{\mathcal{O}^{r}_{C}/C/k}$. Since the whole construction
is algebraic, it can be performed over any complete valued field,
for example, a $p-$adic field.


Let $Q_{d,r,n}(C)$ be the Quot scheme parametrizing rank $r$
coherent sheaf quotients of $\mathcal{O}^{n}_{C}$ of degree $d$. It
is a fine moduli space that comes equipped with a universal exact
sequence over $Q_{d,r,n}(C)\times C$:

$$0\rightarrow \mathcal{K}\rightarrow \mathcal{O}^{n}_{Q_{d,r,n}(C)\times C}
\rightarrow \mathcal{E}\rightarrow 0 $$

\noindent The universal quotient $\mathcal{E}$ is flat over the
$Q_{d,r,n}(C)$ Quot scheme, that is, for each $q\in Q_{d,r,n}(C)$,
$E_{q}:=\mathcal{E}_{|\{q\}\times C}$ is a coherent sheaf over $C$
and
$$h^{0}(E_{q})-h^{1}(E_{q})=d+2\,(1-g),$$ is constant by
Riemann-Roch. That is, it does not depend on $q$. The universal subbundle $\mathcal{K}$ is a locally free sheaf of rank 2, therefore $\varphi: \mathcal{K}\rightarrow \mathcal{O}^{4}$ is a morphism of locally free sheaves. Observe that the sheaf $\mathcal{E}$ is not locally free in the points $(p,t)\in Q_{d,r,n}(C)\times C$, where the rank of the map $\varphi: \mathcal{K}|_{r,t}\rightarrow \mathcal{O}^{4}|_{r,t}$ is 0.

By analogy with $Sym^{d}\,C$, it is natural to define
maps
$$v:Quot^{m}_{\mathcal{O}^{r}_{C}/C/k }\rightarrow J(C),$$ of
Abel-Jacobi type. The geometry of the curve $C$ interacts with the
geometry of $Q_{d,r,n}(C)$ and $J(C)$ via these maps. Note that $J(C)$ is identified with $J^{d}$ by means of a degree $d$ line bundle.

\begin{prop}\label{morphism} \label{morp1} For $d$ sufficiently large and coprime with $r$,
there is a morphism from the Quot scheme $Q_{d,r,n}(C)$ to the
Jacobian of the curve $J^{d}$.
\end{prop}
{\it Proof.} Let $M(r,d)$ be the coarse moduli scheme of stable rank $r$ and degree $d$ vector bundles on $C$ and let $\mathcal{U}$ be the universal bundle over $C\times
M\,(r,d)$. We consider the projective bundle $\rho:
P_{d,r,n}(C)\rightarrow M(r,d)$ whose fiber over a stable bundle
$[F]\in M(r,d)$ is $\mathbb{P}(H^{0}(C,F)^{\oplus n})$. We take
the degree sufficiently large to ensure that the dimension of
$\mathbb{P}(H^{0}(C,F)^{\oplus n})$ is constant. Globalizing, we have
$$P_{d,r,n}(C)=\mathbb{P}(\mathcal{U}^{\oplus n}).$$
Alternatively, $P_{d,r,n}(C)$ may be thought of as a fine moduli
space for $n+1-$pairs  $(F;\phi_{1},\ldots,\phi_{n})$ of a stable
rank $r$, degree $d$ bundle $F$ together with a non-zero $n-$tuple
of holomorphic sections
$\phi=(\phi_{1},\ldots,\phi_{n}):\mathcal{O}^{n}\rightarrow F$
considered projectively. When $\phi$ is surjective, it
defines a point of the Quot scheme $Q_{d,r,n}(C)$,
$$0\rightarrow N\rightarrow \mathcal{O}^{n}\rightarrow E\rightarrow
0$$ where $N=F^{\vee}$. The induced map $\varphi:
Q_{d,r,n}(C)\rightarrow P_{d,r,n}(C)$ is a birational morphism, so
that $Q_{d,r,n}(C)$ and $P_{d,r,n}(C)$ coincide on an open
subscheme and also the universal structures coincide.

From the universal quotient
$$\mathcal{O}^{n}_{Q_{d,r,n}(C)\times C}\rightarrow \mathcal{E}_{Q_{d,r,n}(C)\times C}$$
for all $q\in Q_{d,r,n}(C)$, we have a surjective morphism
$$\mathcal{O}^{n}_{C}\rightarrow E\rightarrow 0.$$

\noindent We now consider the canonical morphism to the Jacobian
of the curve:
$$det: M(r,d)\rightarrow J^{d}.$$ Then the composition of the
morphisms, $\rho=det\circ \rho\circ \varphi$ gives a morphism from
$Q_{d,r,n}$ to the Jacobian $J^{d}$. \cqd

\begin{remark} 
For $L$ a degree $d$ line bundle on $C$, the fiber
$\rho^{-1}([L])$ at $L=\bigwedge^{r}F$ where $[F]\in M(r,d)$, is
isomorphic to $ \mathbb{P}(H^{0}(C,L)^{\oplus\,n})$. In particular,
if
 $r=0$ then  $\rho^{-1}([L])$ corresponds to the variety of higher
rank (n,d)-divisors $E \subseteq \mathcal{O}_{C}(D)^{n}$.

\end{remark}

\begin{remark} The morphism $\rho: Q_{d,r,n}\rightarrow J^{d}$ has good functorial properties, that is, it is compatible with pull-backs, pushforwards and products. In particular, if we consider a cycle class $[\alpha] \in H^{i}(J^{d})$, the pull-back $\rho^{*}([\alpha])$ defines a cycle class in $H^{i}(Q_{d,r,n})$.
\end{remark}

\subsection{Brill-Noether theory on the Quot scheme}


Recalling the notation of \cite{Mar2},  let $R_{C,d}$  be the Quot
scheme compactifying the variety of morphisms $Mor_{d}(C, G(2,4))$,
so that we are fixing the invariants $r$ and $n$ in $Q_{d,r,n}$ to be 2 and 4
respectively. The image of a curve $C$ by $f$ is a geometric curve
in $G(2,4)$ or equivalently a ruled surface in $\mathbb{P}^{3}$. For
each $f:C\rightarrow G(2,4)$ there exists a unique corresponding
quotient $\mathcal{O}^{4}_{C}\rightarrow f^{*}\mathcal{Q}\rightarrow
0$ in $R_{C,d}$, where $\mathcal{Q}$ is the universal quotient over
the Grassmannian. We restrict here to the case of morphisms to the Grassmannian $G(2,4)$,
 because in this case the pull-back $f^{*}\mathcal{Q}$  of the universal quotient bundle, is a rank 2 vector bundle over $C$ and thus there is a well defined invariant called the Segre invariant.

Let us denote by $s$ the Segre invariant $s$ of
the bundle $f^{*}\mathcal{Q}$, which is defined as the minimal
degree of $f^{*}\mathcal{Q}^{\vee}\otimes L$ having a non-zero
section and which satisfies $s\,\equiv d\ (mod\, 2)$ and $2-2g\leq
s\leq g$. In particular, when the minimum value $s$ of $deg(f^{*}\mathcal{Q}^{\vee}\otimes L)$ is achieved for some line subbundle $L$, then $deg\,(L)=\frac{d+s}{2}$. In other words, if $L$ is a line subbundle of $f^{*}\mathcal{Q}$ of degree $\frac{d+s}{2}$, we say it is a maximal line subbundle of $f^{*}\mathcal{Q}$.

Since the universal quotient $\mathcal{E}$ is flat over
$R_{C,d}$, for each $q\in R_{C,d}$,
$E_{q}:=\mathcal{E}|_{\{q\}\times C}$ is a coherent sheaf over $C$ and it is isomorphic to the pull-back $f_{q}^{*}\mathcal{Q}$ of the universal quotient bundle by the corresponding morphism $f_{q}$.
 The variety $\mor_{d}(C, G(2,4))$ of morphisms, sits inside $R_{C,d}$ as the open
subscheme of locally free quotients of $\mathcal{O}^{4}_{C}$.

 Analogously to the case of the Jacobian, we can consider
the following Brill-Noether loci associated with a line bundle $L$
of degree $\frac{d+s}{2}$ on $C$ for a fixed integer $k$:
$$R^{k}_{C,d,s}=\Big\{q\in R_{C,d}|\,h^{0}\,(C,E_{q}^{\vee}\otimes
L)\geq k, \ deg\,L=\frac{d+s}{2}\Big\}=$$

$$\Big\{q\in R_{C,d}|\,h^{1}\,(C,E_{q}\otimes  K_{C}\otimes L^{-1})\geq k, \
deg\,L=\frac{d+s}{2}\Big \}=$$
$$\Big\{q\in R_{C,d}|\,h^{0}\,(C,E_{q}\otimes K_{C} \otimes  L^{-1})\geq k+2g-2-s,
\ deg\,L=\frac{d+s}{2}\Big\}.$$


The subset $R^{k}_{C,d,s}$ has a canonical scheme
structure on $R_{C,d}$, since it can be described as the
degeneration locus of some morphism of vector bundles on
$R_{C,d}$. These sets are analogous to the varieties of special
divisors in the Jacobian of a curve. Note that in the case $k=1$,
this scheme corresponds exactly to the points $q\in R_{C,d}$ such
that $f^{*}_{q}\mathcal{Q}$ has Segre invariant $s$, (this case has been studied in detail in 
\cite{Mar2}). Next proposition shows that when $s$ takes the value $2\,(g-1)$, this defines
a codimension one locus inside $R_{C,d}$ which we will later take
as a polarization for $R_{C,d}$.

\begin{prop} If the Segre invariant $s$ takes the value $2\,(g-1)$, then $R^{1}_{C,d,s}$ is a divisor in $R_{C,d}$.
\end{prop}
{\it Proof.}  By the Krull-Schmidt Theorem, the rank 2 bundle  $f^{*}\mathcal{Q}$ has a unique decomposition into a direct  sum of indecomposable bundles $f^{*}\mathcal{Q} \cong M_{1}\oplus M_{2} $. There are two possibilities:
\begin{enumerate}
\item If $f^{*}\mathcal{Q}$ decomposes into a sum of line bundles, then $M_{1}$ is a line bundle of degree $s$ and $M_{2}$ of degree $d-s$ (resp. $M_{1}$ of degree $d-s$, $M_{2}$ of degree $s$). The condition for the non-emptiness  of the Brill-Noether locus implies $s\leq g$ (see Proposition 3.1 of \cite{Mar1}). As we are asumming $s=2g-2$, this case only occurs when the genus of $C$ is 0, 1 or 2.

\item If $f^{*}\mathcal{Q}$ is indecomposable $M_{1}\cong \O_{C}$ and $M_{2}\cong f^{*}\mathcal{Q}$ (resp. $M_{2}\cong \O_{C}$ and $M_{1}\cong f^{*}\mathcal{Q}$) and since  $2-2g\leq s \leq g$ (see Theorem 2.2 of \S V. 2 of \cite{Har}) and we are assuming that $s=2\,(g-1)$,  the Segre invariant $s=0$. Thus being $s\leq g$,  the  locus $R^{1}_{C,d,s}$ is non-empty as was proved in \cite{Mar1}.
\end{enumerate}
We just need to prove that the stratum $R^{1}_{C,d,s}$ is different from the ambient variety $R_{C,d}$. Then by Theorem 3.2 of \cite{Mar2}, it would follow that $R^{1}_{C,d,s}$ is irreducible and has the right codimension 1.
But the strata $R^{1}_{C,d,s}$ as defined in \cite{Mar2} are dense in $R_{C,d}$, just when the subbundles $M_{1}$ and $M_{2}$ are of maximal degree $\frac{d+s}{2}$ (resp. minimal degree $\frac{d-s}{2}$), in which the dimension of $R_{C,d,s}$ is the maximal one, that is the dimension of $R_{C,d}$. But as we are assuming $s=2g-2$, we have $s\leq g$, and thus by a dimensional argument it follows that $R^{1}_{C,d}$ defines a divisor in $Pic\,(R_{C,d})$.

\cqd

\subsection{Tangent spaces}
Let $0\rightarrow N_{q}\rightarrow \mathcal{O}_{C}^{4}\rightarrow
E_{q}\rightarrow 0$ be the quotient represented by a point $q\in
Q_{d,r,n}(C)$. We consider the tangent space to that point,
$$\mathcal{T}_{q}Q_{d,r,n}(C)\cong Hom\,(N_{q},E_{q})\cong H^{0}(N^{*}_{q}\otimes
E_{q}).$$

If $H^{1}(C,N^{*}_{q}\otimes E_{q})\cong Ext^{1}(N_{q},E_{q})$ is
trivial, then $q$ is a smooth point in $Q_{d,r,n}(C).$ In that
case, we compute using the Riemann-Roch theorem the dimension of
$T_{q}Q_{d,r,n}(C)$ to be $deg\,(N^{*}_{q} \otimes E_{q})+r\cdot
(1-g)$.

\begin{lemma}\label{sing} The singular locus of $R^{k}_{C,d,s}$ contains strictly $R^{k+1}_{C,d,s}$.
\end{lemma}

{\it Proof}. Since the schemes $R^{k}_{C,d,s}$ are determinantal
varieties, locally there exists a morphism $g$ from $R^{k}_{C,d,s}$ to
the variety of matrices $M_{k}(n,m)$ of rank less or equal than $k$,
such that $R^{k}_{C,d,s}$ is the pull-back by $g$ of the variety of
matrices $M_{k}(n,m)$ of rank equal or less than $k$. Then from Proposition II.2 of \cite{ACGH} on  determinantal varieties of matrices, the singular locus of $M_{k}(n,m)$ is exactly equal to $M_{k+1}(n,m)$. Thus  the image of $R^{k}_{C,d,s}$ is a linear space contained in $M_{k}(n,m)$ not meeting $M_{k+1}(n,m)$  which is its singular locus and then the result follows.
\cqd

\begin{lemma}\label{dim}For $d$ sufficiently large depending on $s$,
the expected codimension of $R^{k}_{C,d,s}$ as a determinantal
variety is $(2g-s-2+k)\cdot k$.
\end{lemma}
{\it Proof}. $R^{k}_{C,d,s}$ is exactly the locus of degeneration
of the morphism of vector bundles:
$$\phi:\mathbf{ R}^{1}\pi_{1*}(\mathcal{K}\otimes
\pi^{*}_{2}(L^{-1}\otimes K_{C}))\rightarrow
\mathbf{R}^{1}\pi_{1*}(\mathcal{O}^{n}_{R_{C,d}\times C }\otimes
\pi^{*}_{2}(L^{-1}\otimes K_{C})),$$ where $L$ is a line bundle of
degree $\frac{d+s}{2}$ and $d \equiv s\ mod \ 2$, (see Theorem 3.2
of \cite{Mar2}). Then by the theory of determinantal varieties, a
simple dimension computation gives the result. \cqd

The Segre invariant admits an obvious generalization to higher rank vector bundles. If we consider any Grassmannian $G(r,n)$ of $r-$planes in an $k-$dimensional vector space $V$ over $k$, to any morphism $f:C\rightarrow G(r,n)$ we associate the universal exact sequence over the Grassmannian:

$$0\rightarrow \mathcal{K}\rightarrow \mathcal{O}^{n}_{G}\rightarrow \mathcal{Q}\rightarrow 0.$$
The pull-back $f^{*}\mathcal{Q}$ of the universal quotient bundle by $f$ is now a vector bundle over $C$ of rank $r$. Its Segre invariant $s(r)$ is the maximal degree of $f^{*}\mathcal{Q}^{\vee}\otimes L$ having a non-zero section and satisfying $s\equiv d\, (mod\, r)$ and $2-2g\leq s \leq g$. 

The variety of morphisms $\mor_{d}(C, G(r,k))$ has the expected dimension $nd-r\,(n-r)(g-1)$ if $d>>0$ and moreover in this case, it is irreducible and generically reduced (\cite{Ber}). The Quot scheme compactification $Q_{d,r,n}$ of the variety of morphisms $Mor_{d}(C, G(r,n))$ is endowed with a universal exact sequence on $C\times Q_{d,r,n}$:
$$0\rightarrow \mathcal{N}\rightarrow V\otimes \mathcal{O}_{C\times Q_{d,r,n}}\rightarrow \mathcal{E}\rightarrow 0.$$
We can consider similarly to the case of the Grassmannian $G(2,4)$, the following Brill-Noether loci on $Q_{d,r,n}$:
$$Q^{k}_{d,r,n,s(r)}=\Big\{q\in Q_{d,r,n}|\, \, h^{0}(C,E^{\vee}_{q}\otimes L)\geq k, deg\, L=\frac{d+s(r)}{2}\Big\}.$$

\section{A polarization and Torelli-type result for the variety $ R_{C,d}$}


Let $\mathcal{F}$ be a flat family of coherent sheaves on a relative
smooth projective curve $\pi:\mathcal{C}\rightarrow S$, such that
for each member of the family, the Euler characteristic
$\chi(\mathcal{C}_{s},\mathcal{F}_{s})$ vanishes. We associate to
$\mathcal{F}$ a line bundle $det^{-1}\mathbf{R}\pi_{*}(\mathcal{F})$
(up to isomorphism) equipped with a section $\theta_{\mathcal{F}}$.
We define the theta line bundle on $J$ associated with a line bundle
$L$ of degree $g-1$ on $C$ by applying this construction to the
family $p_{1}^{*}L\otimes \mathcal{P}$ on $C\times J$, where
$\mathcal{P}$ is the Poincare line bundle. Zeroes of the
corresponding theta function $\theta_{L}$ constitute the theta
divisor $\Theta_{L}=\{\xi \in J: h^{0}(L(\xi))>0\},$ which gives
rise to a polarization for the Jacobian $J$. This means that
isomorphism classes of theta line bundles have the form
$det^{-1}\mathcal{S}(L)$, where $\mathcal{S}(L)$ is the
Fourier-Mukai transform of a line bundle $L$ of degree $g-1$
considered as a coherent sheaf on the dual Jacobian $\widehat{J}$
supported on $C$.

Beilinson and Polishchuk gave a proof of the Torelli theorem for
the Jacobian $J$ of a curve in \cite{BP}, based on the observation
that the Fourier-Mukai transform of a line bundle of degree $g-1$
on $C$, is a coherent sheaf (up to shift), supported on the
corresponding theta divisor in $J$. Here we present an analogue of
the Torelli theorem for the variety $R_{C,d}$ of ruled surfaces,
defining a polarization or theta divisor of $R_{C,d}$. Let $C$ be
an algebraic curve of genus $g\geq 2$ and define the map $\nu:
Pic^{\frac{d+s}{2}}(C)\rightarrow Pic^{\frac{g-d-s-1}{2}}(C)$ by
$\nu(L)=K_{C}\otimes L^{-1}$.

In this section we construct an exact functor over the Quot
scheme, and we describe the image of a line bundle on $C$ under this
functor, as a consequence we derive a Torelli
theorem for the Quot scheme.

Let $K_{C}$ be the canonical bundle over $C$ and $\pi_{1}, \pi_{2}$
be the projection maps of $Q_{d,r,n}(C)\times C$ over the first and
second factors respectively. Tensoring the sequence
$$0\rightarrow \mathcal{K}\rightarrow \mathcal{O}^{n}_{Q_{d,r,n}(C)\times C}
\rightarrow \mathcal{E}\rightarrow 0 \ {\rm{over}} \
Q_{d,r,n}(C)\times C\ $$  with the linear bundle $\pi_{2}^{*}(K_{C}\otimes L^{-1})$, yields the exact sequence:

\begin{equation}\label{seq}
0\rightarrow \mathcal{K}\otimes \pi_{2}^{*}(K_{C}\otimes
L^{-1})\rightarrow \mathcal{O}^{n}_{Q_{d,r,n}\times C}\otimes
\pi_{2}^{*}(K_{C}\otimes L^{-1})\rightarrow \mathcal{E}\otimes
\pi_{2}^{*}(K_{C}\otimes L^{-1})\rightarrow 0
\end{equation}

\noindent Here $L$ is a line bundle of fixed degree. The $\pi_{1*}$ direct
image of the above sequence yields the following long exact
sequence on $Q_{d,r,n}(C)$:

$$
0\rightarrow \pi_{1*}(\mathcal{K}\otimes \pi_{2}^{*}(K_{C}\otimes
L^{-1}))\rightarrow \pi_{1*}( \mathcal{O}^{n}_{Q_{d,r,n}}\otimes
\pi_{2}^{*}(K_{C}\otimes L^{-1}))\rightarrow$$ $$\rightarrow
\pi_{1*}(\mathcal{E}\otimes \pi_{2}^{*}(K_{C}\otimes
L^{-1}))\rightarrow \mathbf{R}^{1}\pi_{1*}(\mathcal{K}\otimes
\pi_{2}^{*}(K_{C}\otimes L^{-1}))\rightarrow $$ $$ \rightarrow
\mathbf{R}^{1}\pi_{1*}( \mathcal{O}^{n}_{Q_{d,r,n}}\otimes
\pi_{2}^{*}(K_{C}\otimes L^{-1})) \rightarrow
\mathbf{R}^{1}\pi_{1*}(\mathcal{E}\otimes \pi_{2}^{*}(K_{C}\otimes
L^{-1}))\rightarrow 0.
$$

The universal sheaf $\mathcal{E}$ considered as an object in the
derived category $\mathcal{D}^{b}(Q_{d,r,n}(C)\times C)$ of the product by viewing it as a complex situated in degree 0, determines the exact functor of tensor product:
$$\stackrel{\mathbf{ L}}{\otimes}\mathcal{E}: \mathcal{D}^{b}(Q_{d,r,n}(C)\times C)\rightarrow \mathcal{D}^{b}(Q_{d,r,n}(C)\times C).$$

\noindent The morphism of projection $\pi_{1}$ induces the direct image functor:
$$ \mathbf{R}(\pi_{1*}): \mathcal{D}^{b}(C)\rightarrow \mathcal{D}^{b}(Q_{d,r,n}(C)).$$

\noindent These two functors combined define the exact functor $\phi_{\mathcal{E}}(-)$ with kernel
$\mathcal{E}$ to be the full derived functor 
\begin{eqnarray}
\chi_{\mathcal{E}}(-): \mathcal{D}^{b}(C) \rightarrow
\mathcal{D}^{b}(Q_{d,r,n}(C)) \\ L \mapsto \mathbf{R}\pi_{1*}(\mathcal{E}\stackrel{\mathbf{ L}}{\otimes}
\pi_{2}^{*}(L)),
\end{eqnarray}
between the bounded derived categories of coherent
sheaves over $Q_{d,r,n}(C)$ and $C$ respectively. The same object determines another functor:

\begin{eqnarray}
\psi_{\mathcal{E}}(-): \mathcal{D}^{b}(Q_{d,r,n}(C)) \rightarrow
\mathcal{D}^{b}(C) \\ L \mapsto \mathbf{ R}\pi_{2*}(\mathcal{E}\stackrel{\mathbf{L}}{\otimes}
\pi_{1}^{*}(L)). 
\end{eqnarray}

\begin{rem}The definition of the functor $\chi_{\mathcal{E}}$ uses the
universal quotient sheaf $\mathcal{E}$ which is defined over the
product $Q_{d,r,n}(C)\times C$. However, since the
integral functor $\chi_{\mathcal{E}}$ pushes forward the
corresponding coherent sheaf on the product to $Q_{d,r,n}(C)$,
we can take any $\mathcal{E}\otimes \pi_{1}^{*}N$,
where $N$ is a line bundle over $Quot_{d,r,n}(C)$, to define an
integral exact functor on $\mathcal{D}^{b}(Quot_{d,r,n}(C))$.
\end{rem}

\begin{rem} 
Since $Q_{d,r,n}(C)$ and $C$ have different dimension, the functor
$\chi_{\mathcal{E}}$ cannot be an equivalence between the
corresponding derived categories of coherent sheaves. However it is a representable functor or a functor of Fourier-Mukai type of kernel $\mathcal{E}$.
\end{rem}


\begin{lemma}\label{lem1} The composition of two full and faithful functors is always full and faithful.
\end{lemma}
{\it Proof.} We recall that a functor $F:\mathcal{C}\rightarrow \mathcal{D}$ is said to be fully faithful, if for every pair of objects $X$ and $Y$ in $\mathcal{C}$, the morphism of sets:

$$F: Hom_{\mathcal{C}}(X,Y)\rightarrow Hom_{\mathcal{D}}(F(X),F(Y))$$ is bijective. Now it is clear that if
$F:\mathcal{C}\rightarrow \mathcal{D}$ and $G:\mathcal{D}\rightarrow \mathcal{J}$ are fully faithful functors, the composition $G\circ F: \mathcal{C}\rightarrow \mathcal{J}$ satisfies
$$Hom_{\mathcal{C}}(X,Y)\cong Hom_{\mathcal{D}}(F(X), F(Y))\cong Hom_{E}(GF(X),GF(Y)), $$ and by the definition of composition of functors, the above composite map is $G\circ F$. Thus $G\circ F$ is fully faithful.

\cqd

\begin{prop}\label{compf} The compostion functor $ \psi_{\mathcal{E}}\circ \chi_{\mathcal{E}}: \mathcal{D}^{b}(C)\rightarrow \mathcal{D}^{b}(C)$ is a Fourier-Mukai equivalence.
\end{prop}
{\it Proof.} Given an object $L\in \mathcal{D}^{b}(C)$, we can see by applying the composition $\psi_{\mathcal{E}}\circ \xi_{\mathcal{E}}$ that  $$\psi(\chi(L)):=\mathbf{ R}\pi_{2*}(\mathcal{E}\otimes\pi_{1}^{*}(\mathbf{ R}^{1}\pi_{1*}(\mathcal{E}\otimes \pi_{2}^{*}(L)))).
$$

Now the functor $\pi_{1}^{*}\pi_{1*}: \mathcal{D}(Q_{d,r,n}\times C)\rightarrow \mathcal{D}(Q_{d,r,n}\times C)$ being the composition of two full and faithful functors, the pull-back functor $\pi_{1}^{*}$ and the image direct functor $\pi_{1*}$, is also fully faithful. Then by Orlov's result (Theorem 2.2 of \cite{Orl}), it is represented by an object $\mathcal{L}\in \mathcal{D}(Q\times C)$. This object is constructed by considering a closed embeding $u: Q\hookrightarrow \mathbb{P}^{N}$ in a projective space given by a very ample invertible sheaf on $Q_{d,r,n}$, and constructing an object $\mathcal{L}'\in \mathcal{D}(\mathbb{P}^{N}\times C)$ such that $\mathcal{L}'=(u\times id)_{*}\mathcal{L}$.

\noindent Thus $\pi_{1}^{*}\pi_{1*}(.)=\mathbf{ R}\pi_{2*}(\mathcal{L}\otimes \pi^{*}_{1}(.))$, and by the projection formula applied to the morphism $\pi_{2}$, we have

$$\psi_{\mathcal{E}}\circ \xi_{\mathcal{E}}(L)=\mathbf{ R}\pi_{2*}(\mathcal{E}\otimes\pi_{1}^{*}\mathbf{ R}\pi_{1*}\mathcal{E}\otimes \mathcal{L}\otimes \pi_{1}^{*}(L)).$$

Observe that $\pi_{1}^{*}\mathbf{ R}\pi_{1*}\mathcal{E}$ is a well defined object in $\mathcal{D}(Q\times C)$ with cohomology sheaves $H^{i}(\pi_{1}^{*}\mathbf{ R}\pi_{1*}\mathcal{E})=H^{i}(\mathcal{E})\otimes \mathcal{O}_{Q_{d,r,n}}$.
Finally by recalling the object $\mathcal{N}=\mathcal{E}\otimes\pi_{1}^{*}\mathbf{ R}\pi_{1*}\mathcal{E}\otimes \mathcal{L}$ in $\mathcal{D}(Q\times C)$, we see that the functor $\psi_{\mathcal{E}}\circ \xi_{\mathcal{E}}: \mathcal{D}(C)\rightarrow \mathcal{D}(C)$ is a full and faithful functor represented by $\mathcal{N}$. 
By Serre duality, we have $\mathcal{N}\cong \mathcal{N}^{\vee}\otimes \pi_{2}^{*}K_{C}$. Thus $\psi_{\mathcal{E}}\circ\chi_{\mathcal{E}} $, has the left (right) adjoint functor $(\psi_{\mathcal{E}}\circ \chi_{\mathcal{E}})^{*}$ with kernel $(\mathcal{N}^{\vee})\otimes \pi_{2}^{*}K_{C}$,
defined on objects by $L \mapsto 
\pi_{1*}(\mathcal{N}^{\vee}\otimes\pi_{2}^{*}( K_{C})\otimes \pi_{2}^{*}(L))$ where, $\mathcal{N}^{\vee}:=\underline{R\mathcal{H}om}^{\bullet}(\mathcal{N},
\mathcal{O}_{Q_{d,r,n}(C)\times C}),$ and $K_{C}$ is the canonical sheaf on $C$.

$(\psi_{\mathcal{E}}\circ \chi_{\mathcal{E}})^{*}\circ (\psi_{\mathcal{E}}\circ \chi_{\mathcal{E}})$  is the identity functor, $${\bf id}_{\mathcal{D}(C)}: Ob\ \mathcal{D}(C)\rightarrow Ob\ \mathcal{D}(C) , \ Mor\ \mathcal{D}(C)\rightarrow Mor \ \mathcal{D}(C),$$ 
and the right adjoint functor $(\psi_{\mathcal{E}}\circ \chi_{\mathcal{E}})^{!}_{\mathcal{N}}(.)=K_{Q}[dim\,Q]\otimes \pi_{1*}(\mathcal{N}^{\vee}\otimes (.))$, where $K_{Q}$ is the canonical sheaf on $Q$, gives the identity functor by composing on the right $(\psi_{\mathcal{E}}\circ\chi_{\mathcal{E}})\circ (\psi_{\mathcal{E}}\circ \chi_{\mathcal{E}})^{!}$.

 



\cqd

\begin{lemma} \label{emb}  There is an embedding $i:C\hookrightarrow Q_{d,r,n}$
from the curve $C$ to the Quot scheme $Q_{d,r,n}$.
\end{lemma}
{\it Proof.} 
Consider effective divisors $D_{1},\ldots, D_{r}$ on $C$ and positive integers $d_{1}\ldots, d_{r}$. For each point $p\in C$, we construct the bundle $E_{p}=\bigoplus_{i=1}^{r} \mathcal{O}_{C}(-d_{i}[p]-D_{i})$.Then $E_{p}$ naturally embeds into $\mathbb{C}^{r}\otimes \mathcal{O}_{C}$, $E_{p}\hookrightarrow \mathbb{C}^{r}\otimes \mathcal{O}_{C}$ and $E_{p}\hookrightarrow \mathbb{C}^{n}\otimes \mathcal{O}_{C}$ by picking an $r-$plane in $\mathbb{C}^{n}$. 

Hence dually, we have a morphism $$\mathbb{C}^{n}\otimes \mathcal{O}_{C}\rightarrow E_{p}^{\vee},$$
thus a point in the Quot scheme. Observe that when $\mathbb{C}^{n}\otimes \mathcal{O}_{C}\rightarrow E_{p}^{\vee},$ is not surjective, it produces points in the boundary of the Quot scheme.

\cqd

\begin{rem} 
Note that the morphism given in Lemma \ref{emb} is not canonical and  there can be several embeddings $C\rightarrow Q_{d,r,n}$ depending of the choice of divisors. One can think of the Quot scheme as a space of rank $r$ and degree $d$ matrix of divisors as Bifet in (\cite{Bif}), where it is proved that the fixed point loci for the standard action of the multiplicative group $\mathbb{G}^{r}_{n}$ on the space of divisors are products of symmetric powers of the curve $C$, and sometimes one of this factors will be a copy of $C$.

\end{rem}



From now on, we assume that $r=2$ and $n=4$, that is, $Q_{d,2,4}$ is the Quot scheme $R_{C,d}$, parametrizing rank 2 quotient locally free sheaves of the rank 4 trivial bundle $\mathcal{O}^{4}_{C}$.










\begin{prop}\label{functorF} There is a representable Fourier-Mukai functor
$$ \xymatrix{F: \mathcal{D}(R_{C,d})  \ar[rr]  &  & \mathcal{D}(R_{C,d})  } .$$
\end{prop}
{\it Proof.}
 Let $\mathbf{ L}i^{*}$ be the left derived functor of the embedding $i: C\rightarrow R_{C,d}$ of Lemma \ref{emb}. It is a Fourier-Mukai functor with kernel the structure sheaf $\mathcal{O}_{\Gamma_{i}}$ of the graph $\Gamma_{i}$ of $i: C\hookrightarrow R_{C,d}$. 
 Choose a base point $p\in C$, and take the direct product of $C$ and $R_{C,d}$ over this base point. Then the graph $\Gamma_{i}$ is defined by:
 
 \begin{eqnarray}
\Gamma_{i}: C\to R_{C,d}\times_{p}C\nonumber\\
x \mapsto (i(x),x).
\end{eqnarray}
 
 If we apply the following result due to Mukai (see Proposition 5.10 of \cite{Huy}) to the objects $\mathcal{O}_{\Gamma_{i}}\in \mathcal{D}(R_{C,d}\times C)$ and $\mathcal{E}\in \mathcal{D}(R_{C,d}\times C)$, and the morphisms 
 
 $$
 \xymatrix{ & \ar[ld]_{p12}R_{C,d} \times C \times R_{C,d}  \ar[rd]^{p23}  \ar[d]_{p13} &  \\
 R_{C,d} \times C  & R_{C,d}\times R_{C,d}  &  C\times R_{C,d}  }
 $$
 we have that the composition $F:= \chi_{\mathcal{E}}\circ \mathbf{L}\, i^{*}: \mathcal{D}(R_{C,d})\rightarrow \mathcal{D}(R_{C,d})$  
 
 $$ \xymatrix{ \mathcal{D}(C)  \ar[rr]^{\chi_{\mathcal{E}}}   &  & \mathcal{D}(R_{C,d})
  \\
&  \mathcal{D}(R_{C,d}) \ar[ru]_{\chi_{\mathcal{E}} \circ \mathbf{L}i^{*}} \ar[ul]^{\mathbf{ L}i^{*}}  & } ,$$

\noindent 
 is represented by the object $\mathcal{R}=(\mathbf{R}p_{13*})(p^{*}_{12}\mathcal{O}_{\Gamma_{i}}\otimes p^{*}_{23}\mathcal{E})$.
 


It has the right adjoint functor $F^{!}$ with kernel $\mathcal{R}^{\vee}\otimes K_{R_{C,d}}[dim R_{C,d}]$,
 where $K_{R_{C,d}}$ is the canonical sheaf on $R_{C,d}$. Thus there is a map $id_{R_{C,d}}\rightarrow F^{!}\circ F$. 


\cqd

\begin{rem} Observe that the functor $F$ is not fully faithful since the functor $\chi_{\mathcal{E}}$ is not fully faithful and thus is not an equivalence.
\end{rem}

\begin{rem}
\noindent Note that the existence of the right adjoint functor
 implies the existence of the left adjoint functor 
$F^{*}$ (and reciprocally) by means of the formula:
$$F^{!}=S_{R_{C}}\circ \mathbf{L}F^{*}\circ S^{-1}_{R_{C}},$$
where $S_{R_{C}}=(.)\otimes K_{R_{C,d}}[dim\, R_{C,d}]$, 
is the Serre functor on
$D^{b}(R_{C})$.
\end{rem}

Let us denote by $f$ the composition functor $\psi_{\mathcal{E}}\circ \chi_{\mathcal{E}}: \mathcal{D}^{b}(C)\rightarrow \mathcal{D}^{b}(C)$ of proposition \ref{compf} and by $f^{*}$ its left adjoint functor. 
 As we showed the composition on the left by the adjoint functor gives the identity functor and thus we can talk of the involutive property of $f$ which allows to recover the curve $C$ as we see in the next theorem.
 
 \noindent By Serre duality, the functor 
\begin{eqnarray}
\xymatrix{\nu(-): \mathcal{D}(C) \ar[rr] & &
\mathcal{D}(C)} \\ L \mapsto 
L^{-1}\otimes K_{C},
\end{eqnarray}
defines a Fourier-Mukai involution.

\begin{theorem}\label{tore} For $d>2\,(g-1)$, and $L\in
Pic^{\frac{d+s}{2}}(C)$, the module $\chi_{\mathcal{E}}(\nu(L)))$ in $\mathcal{D}^{b}(R_{C,d})$ is a coherent sheaf $F$ supported by the divisor $R_{C,d,2(g-1)}$ on $R_{C,d}$.
Moreover, 
 the restriction of $F$ to the non-singular
part of this divisor (understood as a polarization for $R_{C,d}$)
is a line bundle and $L$ can be recovered from this line bundle.
\end{theorem}

{\it Proof.} Let $L$ be a line bundle of degree $\frac{d+s}{2}$ on
$C$ $(d\equiv s \ {\rm{mod}}\ 2)$. It is an object in $\mathcal{D}^{b}(C)$ and we can consider 
the module $F:=\chi_{\mathcal{E}}(\nu(L))=	\mathbf{ R}^{1}\pi_{1*}(\mathcal{E}\otimes \pi_{2}^{*}(K_{C}\otimes L^{-1}))$.
We observe that $F$ is a coherent sheaf 
supported on the divisor $R^{1}_{C,d,2(g-1)}$ in $R_{C,d}$, that
is, on the locus of points $q \in R_{C,d}$ such that
$h^{0}\,(C,E^{\vee}_{q}\otimes L)\geq 1$, or dually
$h^{1}(C,E_{q}\otimes K_{C}\otimes L^{-1})\geq 1$. 

Furthermore,
$F$ is the derived pushforward of $\mathcal{E}\otimes\pi^{*}_{2}(K_{C}\otimes L^{-1})$,
so that we can represent it as the cone of the morphism of vector
bundles on $R_{C,d}$, \cite{Mar2}:

\begin{equation}\label{conucleo}\phi: \mathcal{L}_{0}:=
\mathbf{R}^{1}\pi_{1*}(\mathcal{K}\otimes
\pi^{*}_{2}(L^{-1}\otimes K_{C}))\rightarrow
\mathcal{L}_{1}:=\mathbf{R}^{1}\pi_{1*}(\mathcal{O}^{n}_{R_{C,d}\times C
}\otimes \pi^{*}_{2}(L^{-1}\otimes K_{C})),\end{equation} that is,
by the complex $V_{.}=[\mathcal{L}_{0}\rightarrow
\mathcal{L}_{1}]$ in $\mathcal{D}^{b}(R_{C,d})$.

\noindent Since $\chi_{\mathcal{E}}$ is an isomorphism outside of
$R^{1}_{C,d,2(g-1)}$, it is injective and $F=coker \phi$.
Moreover, when $s=2\,(g-1)$, $R^{1}_{C,d,2(g-1)}$ is a divisor and
it is a polarization for $R_{C,d}$. We see that
$$det \phi= det\,\mathbf{R}^{1}\pi_{1*}(\mathcal{E}\otimes
\pi^{*}_{2}(K_{C}\otimes L^{-1}))= $$ $$\ det\left(
\mathbf{R}^{1}\pi_{1*}(\mathcal{K}\otimes \pi^{*}_{2}(L^{-1}\otimes
K_{C}))\right)\otimes det \left(
\mathbf{R}^{1}\pi_{1*}(\mathcal{O}^{n}_{R_{C,d}\times C }\otimes
\pi^{*}_{2}(L^{-1}\otimes K_{C})\right)^{-1}.
$$

The line bundle $det(\chi_{\mathcal{E}}(\nu(L))^{-1}$ on $R_{C,d}$
has a canonical (up to multiplication by a non-zero scalar) global
section 
 $\theta_{\phi}\in H^{0}(R_{C,d},det(\chi_{_{E}}(\nu(L)))^{-1})$ whose zero locus
constitutes a divisor which is supported on $R^{1}_{C,d,2(g-1)}$,
 that is, $det\,F$ is an equation of
$R^{1}_{C,d,2(g-2)}$.

\noindent If we change the resolution $\mathcal{L}_{0}\rightarrow
\mathcal{L}_{1}$ for $F$, the pair $(det\,F,\theta_{\phi})$ gets
replaced by an isomorphic one.

 There is a natural locally closed embedding 
$j: R^{1}_{C,d,2(g-1)}\hookrightarrow R_{C,d}$.

\noindent Let $R_{C,d,s}^{ns}$ be the set $R_{C,d,2(g-1)}\backslash
R^{1}_{C,d,2(g-1)}$. This set corresponds to the non-singular part
of $R_{C,d,s}$
by Lemma \ref{sing}. Let us denote by $j^{ns}$ the induced embedding $R^{ns}_{C,d,s}\hookrightarrow R_{C,d}$.

\noindent The singular locus of $R^{1}_{C,d,2\,(g-1)}$ contains
$R^{2}_{C,d,2\,(g-1)}$ and by Lemma \ref{dim} it has codimension
in $R_{C,d}$ greater than 2 so that $M:=j^{ns*}F$.

Next we prove that $F$ is the pushforward $F=j^{ns}_{*}M$ by $j^{ns}$
of a line bundle $M$ on the non singular part of
$R^{1}_{C,d,2(g-1)}$.

The restriction of $F$ to the non singular part $R^{ns}_{C,d,2(g-1)}$of the divisor
$R_{C,d,2(g-1)}$ is a line bundle from which $F$ can be
recovered by taking the push-forward with respect to the induced
embedding $j^{ns}$. By
the base change theorem for a flat morphism:
$$\mathbf{L}j^{ns*}F\cong \mathbf{R}\pi_{2*}(\mathcal{E}\otimes \pi^{*}_{2}(K_{C}
\otimes L^{-1}))|_{C\times R^{ns}_{C,d,2(g-1)}}.$$
Since $h^{0}(C, E_{q}\otimes K_{C}\otimes L^{-1})=1$ for every $q\in
R^{ns}_{C,d,2(g-1)}$, by applying the base change theorem again, we deduce that
$$rk\,M|_{q}=1 \ {\rm{for \ \ every}} \ q \in R_{C,d,2(g-1)}^{ns}.$$ Since
$R^{ns}_{C,d,2(g-1)}$ is reduced, $M$ is a line bundle on $R^{ns}_{C,d,2(g-1)}$.

We can characterize the set $\mathcal{M}$ of all line bundles on
$R^{ns}_{C,d,2(g-1)}$ in terms of $(R_{C,d},R^{ns}_{C,d,2(g-1)})$. The
set $\mathcal{M}$ has two properties:
\begin{enumerate}
\item  For every $M\in \mathcal{M}$, $M\otimes \nu^{*}M\cong
K_{R^{ns}_{C,d,2(g-1)}}$, where $\nu$ is the map $L\rightarrow
K_{C}\otimes L^{-1}$. \item The class of $M$ generates
the cokernel of the map

\begin{equation}
\label{eq}Pic(R_{C,d})\rightarrow Pic(R^{ns}_{C,d,2(g-1)}).
\end{equation}
\end{enumerate}

Since the Picard variety has the structure of a polarized abelian
variety, the morphism (\ref{eq}) is a homomorphism of abelian
groups induced by the inclusion $j:
R^{ns}_{C,d,s}\hookrightarrow R_{C,d}$. Its cokernel
$coker(Pic\,(R_{C,d})\rightarrow Pic(R^{ns}_{C,d,s}))$ is a
group. Moreover it is isomorphic to $\mathbb{Z}$ for $d$
sufficiently large, (Prop. \ref{morphism}).

Thus to recover the curve $C$ from $(R_{C,d},R_{C,d,2(g-1)})$, by
the involutive property of the integral transform 
$f$, we can recover the curve $C$ by taking a
line bundle $M$ on $\mathcal{M}$, extending it to
$R^{1}_{C,d,2(g-1)}$ by taking the image under the pushforward map $j_{*}$, and then  applying again the integral functor $\psi_{\mathcal{E}}$. Since $F=j^{*}_{ns}M$, if we apply the adjoint functor $f^{*}$ of $f$ to 
$f(\nu(L)):=\psi_{\mathcal{E}}(j^{ns}_{*}(M))$, $(f^{*}\circ f)(\nu(L))$, we recover the line bundle $L$ and thus we get the Torelli result.



 We just need to show that $M \in \mathcal{M}$. First, we apply duality
theory to the projection $\pi_{1}: R_{C,d}\times C\rightarrow
R_{C,d}$ to prove that \begin{equation}
\label{eq1}\mathbf{R}\underline{Hom}(F, \mathcal{O}_{R_{C,d}})\cong
\nu^{*}F[-1],\end{equation} where $\nu: Pic(C)\rightarrow Pic(C)$
corresponds to the involution $L\rightarrow L^{-1}\otimes K_{C}$.
Applying the left derived functor $\mathbf{L}j^{ns*}$ to the isomorphism ($\ref{eq1}$), we
obtain
\begin{equation}\mathbf{R}\underline{Hom}(\mathbf{L}j^{ns*}F,\mathcal{O}_{R^{ns}_{C,d,2(g-1)}})\cong
\nu^{*}\mathbf{L}j^{ns*}F[-1].
\end{equation}
Since $\mathbf{L}j^{ns*}F$ has locally free cohomology sheaves, this
implies that $M^{-1}\cong \nu^{*}\mathbf{L}j^{ns*}F[-1]$. But
$$\mathbf{L}^{-1}j^{ns*}F\cong\mathbf{ L}^{-1}j^{ns*}j^{ns}_{*}M\cong M\otimes
\mathcal{O}_{R^{ns}_{C,d,2(g-1)}}(-R_{C,d,2(g-1)}),$$ and that $\nu^{*}M^{-1}\cong
M(-R_{C,d,2(g-1)})$, which proves condition (1) of $\mathcal{M}$.

\noindent In order to prove the second condition, we consider the
universal quotient $\mathcal{E}|_{\{p\}\times
R^{1}_{C,d,2\,(g-2)}}$ restricted to $\{p\}\times
R^{1}_{C,d,2\,(g-2)}$. The line bundle $M^{-1}$ on
$$R^{ns}_{C,d,2(g-1)}=\{q\in R_{C,d}| h^{0}(C,E^{\vee}_{q}\otimes L)=1, \
L\in Pic^{\frac{d+s}{2}}(C)\}$$ is isomorphic to
$j^{ns}_{*}p_{2*}(\mathcal{O}(R_{C,d,2\,(g-2)}))(-R_{p})$, where
$p_{2}:C\times R_{C,d}\rightarrow R_{C,d}$, $j^{ns}$ is the
embedding $j^{ns}: R^{ns}_{C,d,s}\hookrightarrow R_{C,d}$ and
$R_{p}:=R_{C,d,2\,(g-2)}\cap \{p\} \times R_{C,d}$. Therefore
$M^{-1}\cong \alpha^{*}(\mathcal{O}_{R_{C,d}}(-R_{p}))$ which
generates the cokernel of the map $Pic\,(R_{C,d})\rightarrow
Pic\,(R^{ns}_{C,d,2(g-1)})$. \cqd


\begin{rem}
Observe that we can get the same kind of result by using the Fourier-Mukai transform of Proposition \ref{functorF}.
We see that 
$F(i_{*}(\nu(L))):=\chi_{\mathcal{E}}(i^{*}(i_{*}(\nu(L)))$, by 
 the definition of the functor $F$. 

\end{rem}

\begin{rem}
Note that once we recover the line bundle $L$ supported on $C$, the Torelli result follows from the classical Torelli theorem from which one can recover the curve from the Jacobian $J(C)$. Consider the morphism $\rho$ of Proposition \ref{morphism}. Then the functor $\chi_{\mathcal{E}}$ extends naturally to a functor defined on $\mathcal{D}(J)$,  composing $\chi_{\mathcal{E}}$ by the right with the pull-back functor $\rho^{*}: \mathcal{D}(J)\rightarrow \mathcal{D}(Q)$. 
\end{rem}

\begin{coro} Given two smooth projective curves $C_{1}$ and
$C_{2}$, if there exists an isomorphism \footnote{By isomorphism of Quot schemes as polarized varieties, we mean there is an isomorphism between $R_{C_{1},d}$ and $R_{C_{2},d}$ taking $\theta_{1}$ to $\theta_{2}$} $$f:
(R_{C_{1},d},\theta_{1})\stackrel{\sim}{\rightarrow} (R_{C_{2},d},
\theta_{2})$$ of polarized Quot-schemes, then $C_{1}\simeq C_{2}$.
\end{coro}
{\it Proof.} By Theorem \ref{tore}, the restriction of
$F=\chi_{\mathcal{E}}(\nu_{*}(L))$ to the non-singular part of
$\theta_{i}$ ($L\in Pic^{\frac{d+s}{2}}$), is a line bundle
$j^{ns}_{*}M$ and $L$ can be recovered from this line bundle since
$M:=j^{ns*}F$. Therefore
$$f|_{supp(j_{1*}^{ns})}:
supp(j_{1*}^{ns}M_{1})\stackrel{\sim}{\rightarrow}
supp(j_{2*}^{ns}M_{2}),$$ and $C_{1}\cong C_{2}$, where
$j_{1}:\theta_{1}^{ns}\hookrightarrow R_{C_{1},d}$, and
$j_{2}:\theta_{2}^{ns}\hookrightarrow R_{C_{2},d}$. \cqd

\subsubsection*{Acknowledgments}
The idea to use a Fourier-Mukai functor applied to this kind of
problem, was inspired by the work of  A. Polishchuk and A. Beilinson,
\cite{BP}.


\end{document}